\documentstyle[psfig,intlim,12pt]{amsart}
\newcommand{\C}{\Bbb C}
\newcommand{\R}{\Bbb R}
\newcommand{\Z}{\Bbb Z}

\newtheorem{th}{Theorem}
\newtheorem{lemma}{Lemma}
\newtheorem{thm}{Theorem}[section]
\newtheorem{lem}[thm]{Lemma}

\theoremstyle{definition}

\theoremstyle{remark}

\newtheorem{ex}[thm]{Example}
\newtheorem{rmk}{Remark}
\newtheorem{ack}{Acknowledgment}

\newcommand{\degree}{\operatorname{deg}}
\newcommand{\conj}{\operatorname{conj}}
\newcommand{\dd}{\partial}

\begin{document}
\title[Topology of curves of degree six on cubics]
{Topological arrangement of curves of degree 6
on cubic surfaces in \protect$\R P^3$}
\author{G.Mikhalkin}
\thanks{Research at MSRI is supported in part by NSF grant
DMS-9022140}



\address{Department of Mathematics\\
University of Toronto\\
Toronto, Ont. M5S 3G3\\ Canada}
\email{mihalkin@@math.toronto.edu, grisha@@msri.org}
\curraddr{Mathematical Science Research Institute\\
1000 Centennial Dr.\\
Berkeley, CA 94720, USA}

\maketitle

\begin{abstract}
A quadric in $\R P^3$ cuts a curve of degree 6
on a cubic surface in $\R P^3$.
The papers classifies the nonsingular curves
cut in this way on non-singular cubic surfaces
up to homeomorphism.

Two issues new in the study related to
the first part of the 16th Hilbert problem appear
in this classification.
One is the distribution of the components of the curve
between the components of the non-connected cubic
surface which turns out to depend on the patterns of
arrangements (see Theorem 1).
The other is presence of positive genus among
the components of the complement and genus-related
restrictions (see Theorems 3 and 4).
\end{abstract}

\section{Introduction}
\label{intro}

Let $\R B$ be a real algebraic cubic surface in $\R P^3$ and
let $\R Q$ be a real algebraic quadric surface in $\R P^3$.
This just means that for some homogeneous polynomials
$f$, $\degree(f)=3$, and $g$, $\degree(g)=2$,
 surfaces $\R B$ and $\R Q$ are given by the following equations
\begin{gather*}
\R B=\{[x_0:x_1:x_2:x_3]\in\R P^3\ |\ f(x_0,x_1,x_2,x_3)=0\}\\
\R Q=\{[x_0:x_1:x_2:x_3]\in\R P^3\ |\ g(x_0,x_1,x_2,x_3)=0\}.
\end{gather*}
Let $\C B$ and $\C Q$ be the complexifications of surfaces $\R B$ and $\R Q$
\begin{gather*}
\C B=\{[x_0:x_1:x_2:x_3]\in\C P^3\ |\ f(x_0,x_1,x_2,x_3)=0\}\\
\C Q=\{[x_0:x_1:x_2:x_3]\in\C P^3\ |\ g(x_0,x_1,x_2,x_3)=0\}.
\end{gather*}
We suppose that $\C B$ and $\C Q$ are nonsingular and
intersect transversely along a curve $\C A$.
The intersection $\C A=\C B\cap\C Q$ is a smooth curve in $\C P^3$,
its real part $\R A$ is $\R B\cap\R Q$.

Consider now the diffeomorphism type of $\C B$, $\C Q$, $\C A$.
It is well-known that $\C Q$ is diffeomorphic to
$\C P^1\times\C P^1$,
$\C B$ is diffeomorphic to
$\C P^2\#6\overline{\C P}^2$ and we can construct $\C B$
from $\C P^2$ by blowing up 6 points not sitting on the same conic.
The curve $\C A$ is the proper transform of a plane curve of degree 6
with 6 ordinary double points in the blowup points
so the genus of $\C A$ is 4.

Consider now the diffeomorphism type of $\R B$, $\R Q$, $\R A$.
It is well-known that a quadric $\R Q$ is a hyperboloid
(i.e. diffeomorphic to $\R P^1\times\R P^1$), an ellipsoid
(i.e. diffeomorphic to $S^2$) or empty.
It is also well-known that a cubic $\C B$ is diffeomorphic
to $\#(2n+1)\R P^2\approx\R P^2\#nT^2$, where $n\le 3$ or $\R P^2\sqcup S^2$.
Since $\R A$ is smooth, every component of $\R A$ is diffeomorphic to $S^1$.
Note that the embedding of this $S^1$ into $\R B$ is two-sided,
the sides of $\R A$ are defined by the sign of polynomial $g$
(since $\deg g$ is even).

Since the genus of $\C A$ is 4,
the Harnack inequality implies that the number of components of $\R A$
is at most 5 ($\C A/\conj$ is a connected surface and, therefore,
the number of components of its boundary is at most
$2-\chi(\C A/\conj)=2-\frac{1}{2}\chi(\C A)=5$).

The topological classification of pairs $(\R Q, \R A)$
was done by Gudkov \cite{G} in the case where $\R Q$ is a hyperboloid
and by Gudkov and Shustin \cite{GSh} in the case where $\R Q$ is an ellipsoid
(see \cite{Mi2} for the classification of flexible curves).

This paper gives the topological classification
of pairs $(\R B,\R A)$.
Two new phenomena appear in this classification ---
the arrangement of $\R A$ with respect to the components of $\R B$
in the case when $\R B$ is not connected
and the arrangement of the components of $\R A$ with respect to the handles
of $\R B$ in the case when $\R B$ is diffeomorphic to a projective plane
with handles.
These possibilities for arrangement turn out to be
controlled by the degree of the curve.

The topological type of a pair $(\R B,\R A)$ is also called
the type of {\em topological arrangement} of $\R A$ in $\R B$ or
the {\em real scheme} of $\R A$ in $\R B$.

\begin{ack}
A part of the work was done during my stay
at the Institute for Advanced Study
(supported by the NSF grant DMS 9304580)
and at the Max-Planck-Institut.
\end{ack}

\section{Curves on cubic surfaces of positive Euler characteristic}
There are two diffeomorphism types of cubic surfaces
of positive Euler characteristic in $\R P^3$:
$$\R P^2\sqcup S^2\ \ \text{and}\ \ \R P^2.$$
We consider them separately since if $\R B$ is one of these surfaces
than any component of $\R A$ bounds a disk in $\R B$ and even the
notations for the topological type of $(\R B,\R A)$ are simpler in this case.

\subsection{Notations}
Recall the system of codes for topological arrangements of
a collection of two-sided circles (or {\em ovals}) in $\R P^2$
introduced by Viro \cite{V}.
Every oval separates $\R P^2$ into a disk (caller the {\em interior} of the oval)
and a M\"obius band.
Two ovals are called {\em disjoint} if their interiors are disjont.
The system of codes is inductive, $n$ disjoint ovals are encoded by $<n>$
and if we add a new oval containing a collection $<A>$ in its interior
then the code of the new collection is $<1\!<\! A\!>>$.

To make the same system suitable for encoding of topological arrangements
of embedded curves in $S^2$ it suffices to fix a point in $S^2$ in
the complement of the curve --- then the interior of an oval
is defined as the component of the complement of the oval
not containing the fixed point.
Such a code determines uniquely a topological type of curve in $S^2$,
but because of ambiguity in fixing a point there is more than one
code corresponding to the same topological type.
We choose codes minimizing the number of ovals in the interior of other ovals.
To encode types of topological arrangement of $\R A$ into a cubic surface
diffeomorphic to $\R P^2\sqcup S^2$ we combine systems for $\R P^2$ and $S^2$.

Parameters $\alpha$ and $\beta$ in the theorems are nonnegative integer numbers.
\subsection{Classification}
\begin{th}
\label{hyp}
Topological arrangement of a nonsingular algebraic curve
of degree 6 on a nonsingular cubic surface diffeomorphic to $\R P^2\sqcup S^2$
is of one of the following types:
\begin{itemize}
\begin{description}
\item[a] $<\alpha\sqcup1\!<\! 1\!>>_{\R P^2}\sqcup<\emptyset>_{S^2},\hfill
\alpha\le 3,$
\item[b] $<1\!<\!\alpha\!>>_{\R P^2}\sqcup<\emptyset>_{S^2},\hfill
2\le\alpha\le 4,$
\item[c] $<\alpha>_{\R P^2}\sqcup<\beta>_{S^2},\hfill
\alpha+\beta\le 5,$
\item[d] $<1\!<\!1\!<\!1\!>>>_{\R P^2}\sqcup<\emptyset>_{S^2},$
\item[e] $<1\!<\!1\!>>_{\R P^2}\sqcup<1>_{S^2},$
\item[f] $<\emptyset>_{\R P^2}\sqcup<1\sqcup1\!<\! 1\!>>_{S^2}.$
\end{description}
\end{itemize}
Each of the 31 types listed above is realizable
by a nonsingular algebraic curve of degree 6.
\end{th}

\begin{th}
\label{rp2}
Topological arrangement of a nonsingular algebraic curve of degree 6
on a nonsingular cubic surface diffeomorphic to $\R P^2$
is of one of the following types:
\begin{itemize}
\begin{description}
\item[a] $<\alpha\sqcup 1\!<\!\beta\!>>,\hfill \alpha+\beta\le 4,$
\item[b] $<1\!<\!1\!<\!1\!>>>,$
\item[c] $<\emptyset>$.
\end{description}
\end{itemize}
Each of the 17 types listed above is realizable
by a nonsingular algebraic curve of degree 6.
\end{th}

\section{Curves on cubic surfaces of negative Euler characteristic}
There are three diffeomorphism types of cubic surfaces
of negative Euler characteristic in $\R P^3$:
$$\R P^2\#3T^2,\ \ \R P^2\#2T^2\ \ \text{and}\ \ \R P^2\#T^2.$$

\subsection{Notations}
\label{notneg}
We define
$$B_+=\{x\in\R B\ |\ g(x)\ge 0\},\
B_-=\{x\in\R B\ |\ g(x)\le 0\}.$$
Then $\dd B_+=\dd B_-=\R A$.
Note that at least one of $B_+$, $B_-$ is non-orientable since
$$\chi (\R B)\equiv 1\pmod{2}.$$
Changing the sign of $g$ we may assume that
\begin{gather*}
\chi(B_+)\equiv b_0(\R A)\pmod{2}\\
\chi(B_-)\equiv 1+b_0(\R A)\pmod{2},
\end{gather*}
where $b_0$ stands for the number of components (the 0-dimensional Betti number).
Then $B_-$ is always nonorientable.

Note that the topological type of $(\R B,\R A)$ is determined by
the topological types of $B_+$ and $B_-$ if
the number of non-disk components is not more than one for
both $B_+$ and $B_-$ (since then there is only one way to glue
$B_+$ and $B_-$ along the boundary to obtain a connected $\R B$).
In this case we code the topology of $\R A$ on $\R B$ with
$$<j\sqcup F,n\sqcup G>$$
where $j$ is the number of disks in $B_+$, $n$ is
the number of disks in $B_-$,
$F$ is the topological type
of the only non-disk component of $B_+$ (if any)
and $G$ is the only non-disk component of $B_-$.
We denote the topological type of the connected
orientable surface of genus $g$ with $gT^2$ and
that of the non-orientable one with $g\R P^2$.
Notation $F_k$ stands for the result of puncturing
of a closed surface $F$ $k$ times.

Lemma \ref{arnold} implies that if $B_+$ or $B_-$
has more than one non-disk components then the
topological type of $(\R B,\R A)$ is still determined
by the topological types of $B_+$, $B_-$ and $\R B$.
If $\R B\approx\R P^2\# nT^2$, $n=1,2,3$ then
the code of $(\R B,\R A)$ is
$$<1\sqcup nT^2_2,S^2_2\sqcup\R P^2_1>$$
if the core of the annulus is null-homologous and
$$<(n-1)T^2_3,S^2_2\sqcup\R P^2_1>$$ if not.

\subsection{Classification}
\begin{th}
\label{7rp2}
Topological arrangement of a nonsingular algebraic curve
of degree 6 on a nonsingular cubic surface diffeomorphic to $\R P^2\#3T^2$
is of one of the following types:
\vspace{10pt}\newline
{\bf 1.} $<\alpha\sqcup S^2_{\beta+\gamma},
\beta\sqcup(9-2\gamma)\R P^2_{\alpha+\gamma}>$,\\
\begin{tabular}{|r|r|r|r|r|}
\hline
\ & $\beta=0$ & $\beta=1$ & $\beta=2$ & $\beta=3$\\
\hline
$\alpha=0$ & $\gamma=1,2,3,4$ & $\gamma=1,2,3,4$ &
$\gamma=1,2,3$ & $\gamma=2$\\
\hline
$\alpha=1$ & $\gamma=1,2,3,4$ & $\gamma=1,2,3$ &
$\gamma=2$ & \ \\
\hline
$\alpha=2$ & $\gamma=1,2,3$ & $\gamma=1,2$ & \ & \ \\
\hline
$\alpha=3$ & $\gamma=1,2$ & $\gamma=1$ & \ & \ \\
\hline
\end{tabular}\vspace{10pt}\\
{\bf 2.} $<\alpha\sqcup T^2_{\beta+\gamma},
\beta\sqcup(7-2\gamma)\R P^2_{\alpha+\gamma}>,$\\
\begin{tabular}{|r|r|r|r|r|r|}
\hline
\ & $\beta=0$ & $\beta=1$ & $\beta=2$ & $\beta=3$ & $\beta=4$ \\
\hline
$\alpha=0$ & $\gamma=1,2,3$ & $\gamma=1,2,3$ & $\gamma=1,2,3$ &
$\gamma=1,2$ & $\gamma=1$\\
\hline
$\alpha=1$ & $\gamma=1,2,3$ & $\gamma=1,2,3$ & $\gamma=1,2$ & \ & \ \\
\hline
$\alpha=2$ & $\gamma=1,2,3$ & $\gamma=1,2$ & $\gamma=1$ & \ & \ \\
\hline
$\alpha=3$ & $\gamma=1,2$ & \ & \ & \ & \ \\
\hline
$\alpha=4$ & $\gamma=1$ & \ & \ & \ & \ \\
\hline
\end{tabular}\vspace{10pt}\\
{\bf 3.} $<\alpha\sqcup 2T^2_{\beta+\gamma},
\beta\sqcup(5-2\gamma)\R P^2_{\alpha+\gamma}>$,\\
\begin{tabular}{|r|r|r|r|r|}
\hline
\ & $\beta=0$ & $\beta=1$ & $\beta=2$ & $\beta=3$ \\
\hline
$\alpha=0$ & $\gamma=1,2$ & $\gamma=1,2$ & $\gamma=1,2$ &
$\gamma=1,2$ \\
\hline
$\alpha=1$ & $\gamma=1,2$ & $\gamma=1,2$ & $\gamma=1,2$ &
$\gamma=1$ \\
\hline
$\alpha=2$ & $\gamma=1,2$ & $\gamma=1,2$ & \ & \ \\
\hline
$\alpha=3$ & $\gamma=1,2$ & $\gamma=1$ & \ & \ \\
\hline
\end{tabular}\vspace{10pt}\\
{\bf 4.} $<\alpha\sqcup 3T^2_{\beta+1},
\beta\sqcup\R P^2_{\alpha+1}>$,\\
\begin{tabular}{|r|r|r|r|r|}
\hline
$\beta=0$ & $\beta=1$ & $\beta=2$ & $\beta=3$ & $\beta=4$\\
\hline
$\alpha=0,1,2,3$ & $\alpha=0$ & $\alpha=0$ & $\alpha=0$ &
$\alpha=0$ \\
\hline
\end{tabular}\vspace{10pt}\\
{\bf 5.} $<\alpha\sqcup 2\R P^2_{\beta+\gamma},
\beta\sqcup (7-2\gamma)\R P^2_{\alpha+\gamma}>$,\\
\begin{tabular}{|r|r|r|r|r|r|}
\hline
\ & $\beta=0$ & $\beta=1$ & $\beta=2$ & $\beta=3$ & $\beta=4$\\
\hline
$\alpha=0$ & $\gamma=1,2,3$ & $\gamma=1,2,3$ & $\gamma=1,2,3$ &
$\gamma=1,2$ & $\gamma=1$\\
\hline
$\alpha=1$ & $\gamma=1,2,3$ & $\gamma=1,2,3$ & $\gamma=1,2$ &
$\gamma=1$ & \ \\
\hline
$\alpha=2$ & $\gamma=1,2,3$ & $\gamma=1,2$ & \ & \ & \ \\
\hline
$\alpha=3$ & $\gamma=1,2$ & $\gamma=1$ & \ & \ & \ \\
\hline
$\alpha=4$ & $\gamma=1$ & \ & \ & \ & \ \\
\hline
\end{tabular}\vspace{10pt}\\
{\bf 6.} $<\alpha\sqcup 4\R P^2_{\beta+\gamma},
\beta\sqcup (5-2\gamma)\R P^2_{\alpha+\gamma}>$,\\
\begin{tabular}{|r|r|r|r|r|r|}
\hline
\ & $\beta=0$ & $\beta=1$ & $\beta=2$ & $\beta=3$ & $\beta=4$\\
\hline
$\alpha=0$ & $\gamma=1,2$ & $\gamma=1,2$ & $\gamma=1,2$ &
$\gamma=1,2$ & $\gamma=1$\\
\hline
$\alpha=1$ & $\gamma=1,2$ & $\gamma=1,2$ & $\gamma=1,2$ &
$\gamma=1$ & \ \\
\hline
$\alpha=2$ & $\gamma=1,2$ & $\gamma=1,2$ & $\gamma=1$ & \ & \ \\
\hline
$\alpha=3$ & $\gamma=1,2$ & $\gamma=1$ & \ & \ & \ \\
\hline
$\alpha=4$ & $\gamma=1$ & \ & \ & \ & \ \\
\hline
\end{tabular}\vspace{10pt}\\
{\bf 7.} $<\alpha\sqcup 6\R P^2_{\beta+1},
\beta\sqcup \R P^2_{\alpha+1}>$,\\
\begin{tabular}{|r|r|r|r|r|}
\hline
$\beta=0$ & $\beta=1$ & $\beta=2$ & $\beta=3$ & $\beta=4$\\
\hline
$\alpha=0,1,2,3,4$ & $\alpha=0,1$ & $\alpha=0,1$ &
$\alpha=0,1$ & $\alpha=0$\\
\hline
\end{tabular}\vspace{10pt}\\
{\bf 8.}
$<1\sqcup 3T^2_2,S^2_2\sqcup\R P^2_1>$,
$<2T^2_3,S^2_2\sqcup\R P^2_1>$\vspace{3pt}\\
{\bf 9.} $<\emptyset,7\R P^2>$.\vspace{3pt}\\
Each of the 157 types listed above is realizable
by a nonsingular algebraic curve of degree 6.
\end{th}

\begin{th}
\label{5rp2}
Topological arrangement of a nonsingular algebraic curve
of degree 6 on a nonsingular cubic surface diffeomorphic
to $\R P^2\#2T^2$ is of one of the following types:
\vspace{3pt}\\
{\bf 1.} $<\alpha\sqcup S^2_{\beta+\gamma},
\beta\sqcup (7-2\gamma)\R P^2_{\alpha+\gamma}>$,\hfill
$\alpha+\beta+\gamma\le 5$, $\gamma=1,2,3$,\vspace{3pt}\\
{\bf 2.} $<\alpha\sqcup T^2_{\beta+\gamma},
\beta\sqcup (5-2\gamma)\R P^2_{\alpha+\gamma}>$,\hfill
$\alpha+\beta+\gamma\le 5$, $\gamma=1,2$,\vspace{3pt}\\
{\bf 3.} $<\alpha\sqcup 2T^2_{\beta+1},
\beta\sqcup\R P^2_{\alpha+1}>$,\hfill
$\alpha+\beta\le 4$, $\alpha\neq 4$,\vspace{3pt}\\
{\bf 4.} $<\alpha\sqcup 2\R P^2_{\beta+\gamma},
\beta\sqcup (5-2\gamma)\R P^2_{\alpha+\gamma}>$,\hfill
$\alpha+\beta+\gamma\le 5$, $\gamma=1,2$,\vspace{3pt}\\
{\bf 5.} $<\alpha\sqcup 4\R P^2_{\beta+1},
\beta\sqcup\R P^2_{\alpha+1}>$, \hfill
$\alpha+\beta\le 4$,\vspace{3pt}\\
{\bf 6.}
$<1\sqcup 2T^2_2,S^2_2\sqcup\R P^2_1>$,
$<T^2_3,S^2_2\sqcup\R P^2_1>$\vspace{3pt}\\
{\bf 7.} $<\emptyset,5\R P^2>$.\vspace{3pt}\\
Each of the 113 types listed above is realizable by
a nonsingular algebraic curve of degree 6.
\end{th}

\begin{th}
\label{3rp2}
Topological arrangement of a nonsingular algebraic curve
of degree 6 on a nonsingular cubic surface diffeomorphic to $\R P^2\#T^2$
is of one of the following types:
\vspace{3pt}\\
{\bf 1.} $<\alpha\sqcup S^2_{\beta+\gamma},
\beta\sqcup (5-2\gamma)\R P^2_{\alpha+\gamma}>$,\hfill
$\alpha+\beta+\gamma\le 5$, $\gamma=1,2$,\vspace{3pt}\\
{\bf 2.} $<\alpha\sqcup T^2_{\beta+1},
\beta\sqcup\R P^2_{\alpha+1}>$,\hfill
$\alpha+\beta\le 4$,\vspace{3pt}\\
{\bf 3.} $<\alpha\sqcup 2\R P^2_{\beta+1},
\beta\sqcup\R P^2_{\alpha+1}>$,\hfill
$\alpha+\beta\le 4$,\vspace{3pt}\\
{\bf 4.}
$<1\sqcup T^2_2,S^2_2\sqcup\R P^2_1>$,
$<S^2_3,S^2_2\sqcup\R P^2_1>$\vspace{3pt}\\
{\bf 5.} $<\emptyset,3\R P^2>$.\vspace{3pt}\\
Each of the 58 types listed above is realizable by
a nonsingular algebraic curve of degree 6.
\end{th}

\section{Auxiliary lemma}
\begin{lemma}[cf. \protect\cite{A}]
\label{arnold}
$B_+$ contains not more than one component
of non-positive Euler characteristic.
If $B_-$ containes more then one component
of non-positive Euler characteristic
then $B_-\approx S^2_2\sqcup\R P^2_1$
and the cores of $S^2_2$ and $\R P^2_1$ of $B_-$
form a cycle dual to $w_1(\R B)$.
\end{lemma}
\begin{rmk}
Of course, if $\R A$ is empty and $\R B$ is not $S^2\sqcup\R P^2$
then all components of positive Euler characteristic
are disks.
\end{rmk}
\begin{pf} Consider the double
covering $\C X$ of $\C B$ branched along $\C A$.  The complex
conjugation of $\C B$ can be lifted in two different ways to an
involution on $\C Y$ (which differ by the covering automorphism).
The fixed point set $Y_+$ and $Y_-$ of these two liftings $c_+$ and
$c_-$ are double branched coverings of $B_+$ and $B_-$.  Note that
$Y_+$ and $Y_-$ are orientable (even though $B_-$ is never orientable
and $B_+$ may be non-orientable) since $w_1(\R B)$ is induced by the
plane section in $\R P^3$.

The homology vector space $H_2(\C Y;\R)$ splits into
the orthogonal direct sum of the subspaces
$H_2^+(\C Y;\R)$ and $H_2^-(\C Y;\R)$ invariant under actions
of both $c_+$ and $c_-$, where the intersection form of $\C Y$
is positive definite on $H_2^+(\C Y;\R)$ and negative definite on
$H_2^-(\C Y;\R)$.
An easy computation shows that
$$\operatorname{dim}H_2^+(\C Y;\R)=3$$
($\C Y$ is a K3-surface).
Furthermore, $H_2^+(\C Y;\R)$ splits
under the action of $c_+$ and $c_-$ into ($-1$)-eigenspaces
$E^{-1}_{c_+}$ and $E^{-1}_{c_-}$ and ($+1$)-eigenspaces
$E^{+1}_{c_+}$ and $E^{+1}_{c_-}$.
Since $c_+\circ c_-$ is the covering authomorphism
\begin{gather*}
\operatorname{dim}(E^{+1}_{c_+}\cap E^{-1}_{c_-})=1\\
\operatorname{dim}(E^{-1}_{c_+}\cap E^{+1}_{c_-})=1.
\end{gather*}
Note that $c_-$ reverses orientation of $Y_+$.
Therefore if the linear combination of the homology classes realized
by components of $Y_+$ is in $H_2^+(\C Y;\R)$ then it
is in $(E^{+1}_{c_+}\cap E^{-1}_{c_-})$.
Similarily, if the linear combination of the homology classes realized
by components of $Y_-$ is in $H_2^+(\C Y;\R)$ then it
is in $(E^{-1}_{c_+}\cap E^{+1}_{c_-})$.

Note that the self-intersection number of a component of $Y_+$ ($Y_-$)
in $\C Y$ is equal to $(-1)$ times its Euler characteristic since
the multiplication by $i$ provides an orientation-reversing
isomorphism between the tangent and the normal bundles of $Y_+$ ($Y_-$)
in $\C Y$.
Therefore, the number of the components of $Y_+$ ($Y_-$) with
non-positive Euler characteristic is not greater then
$\operatorname{dim}(E^{+1}_{c_+}\cap E^{-1}_{c_-})=1$
($\operatorname{dim}(E^{-1}_{c_+}\cap E^{+1}_{c_-})=1$)
plus the number
of linear relations on the $\R$-homology
classes of these components.

Since $H_1(\C Y)=0$ the universal coefficients theorem implies that
the number of linear relations on the
$\R$-homology classes of the components of $Y_+$ ($Y_-$) with
non-positive Euler characteristic is not more than the number of
corresponding $\Z_2$-relations.
The Smith exact sequence (see e.g.
appendix in \cite{W}) assures that the only possible
non-trivial linear relation
on the $Z_2$-homology classes of $Y_+$ ($Y_-$)
is the sum of all of them.
Therefore the maximal number of the $\R$-relations is 1 and
it can be reached only if all the components of $Y_+$ ($Y_-$)
are of non-positive Euler characteristic.
Furthermore, if at least one of them is of positive Euler
characteristic then they are linearily independent over $\R$
since the only possible $\R$-relation (the sum of them taken with
non-zero coefficients) intersects that component non-trivially.

Therefore, if there is more than one component of $B_+$ ($B_-$)
of non-positive Euler characteristic then $B_+$ ($B_-$) consists
of two components of zero Euler characteristic.
In addition $Y_+$ ($Y_-$) is $\Z_2$-homologous to zero in $\C Y$.
The latter implies that $\R A$ is $\Z_2$-homologous to zero in
$\C A$.
Indeed, if there exists a loop in $\C A$ which
intersects $\R A$ once then a disk bounded by
this loop in $\C B$ (recall that $\pi_1(\C B)=0$)
is a closed surface which intersects $Y_+$ ($Y_-$) in odd
number of points.
This implies in turn that the number of components of $\R A$
is of the opposite pairity with the genus of $\C A$ since
$\R A$ separates $\C A$ into two diffeomorphic complex
conjugate orientable surfaces.

There are two types of connected surfaces of zero Euler
characteristic with boundary --- annuli and M\"obius bands.
The above discussion implies that the only possibility
for the surface $B_+$ ($B_-$) to have more than one component
of non-positive Euler characteristic is to consist of
an annulus and a M\"obius band.
By our convention (see \ref{notneg})
this surface is not $B_+$.

If the Euler characteristic of $\R B$ is positive
($\R B\approx S^2\sqcup\R P^2$ or $\R B\approx\R P^2$)
then the proof of the Lemma is finished since $w_1(\R B)$ is
the only non-trivial element in $H_1(\R B;\Z_2)$.
To finish the proof for $\R B$ of negative Euler characteristic
we consider the birational equivalence $\beta:\C B\to\C P^2$
which is the blowup in 6 points.
The image of $\R A$ under $\beta$ is a curve $\R C$ of degree 6 in $\R P^2$.
which has singularities of multiplicity 2 at the points of blowup.
We have to prove that the image of the union of the cores
of the annulus and the M\"obius band
is dual to $w_1(R P^2)$ and contains all points of blowup.

We deduce it from the Bezout theorem.
The cores are disjoint so if the image of the cores is not
dual to $w_1(\R P^2)$ then the image of the each core bounds
a disk in $\R P^2$.
Pick a point in each of these disks which does not belong
to the image of $B_-$ and draw a line through these points.
This line intersects $\R C$ at least in 8 points which
contradicts to the Bezout theorem.
Therefore, the image of exactly one of the cores bounds a disk.
If there is a blowup point on $\R P^2$ which does not belong to
the image of either of the cores then we draw a line through
that blowup point and a point on the disk which does not belong
to $B_-$.
This line intersects $\R C$ at least in 6 points other than the
blowup point (which is of multiplicity 2).
\end{pf}

\section{Proof of Theorem 1}
The proof of Theorem 1 (as well as the proofs of Theorems 2,
3, 4 and 5) consists
of two parts --- the part "Restrisctions" where we prove that the
topological arrangements not listed in the theorem are not realizable
and the part "Constructions" where we construct the curves with the
listed topological arrangements.
\subsection{Restrictions}
Recall (see section 1) that by the Harnack inequality
the maximal number of components of $\R A$ is 5.
Lemma \ref{arnold} implies that
unless $(\R B,\R A)$ is of type d of Theorem 1
there is no more than one
component of $\R A$ which does not bound a disk in $\R B$
disjoint from the other components of $\R A$.
If this component belongs to the sphere
component of $\R B$ then Lemma \ref{arnold}
implies that $(\R B,\R A)$ is of type $<\emptyset>_{\R P^2}\sqcup
<1\sqcup 1\!<\!1>>_{S^2}$.
If this component belongs to the plane component of $\R B$
then Lemma \ref{arnold} implies that there is no more than one
component of $\R A$ on the sphere component of $\R B$.

Thus we have only to prove that the following arrangements
do not appear
\begin{gather*}
<\alpha\sqcup 1\!<\!\beta\!>>_{\R P^2}\sqcup<\emptyset>_{S^2},\
\text{if $\alpha>0$, $\beta>1$,}\\
<\alpha\sqcup 1\!<\!\beta\!>>_{\R P^2}\sqcup<1>_{S^2},\
\text{if $\beta>0$, unless $\alpha=0$, $\beta=1$.}
\end{gather*}

\begin{lem} $\R B$ is of type I$rel$,
i.e. $\R B$ is dual to the second Stiefel-Whitney class $w_2(\C B)$
of $\C B$.
\end{lem}

\begin{pf} The sphere component of $\R B$
bounds a ball in $\R P^3$.
Let $x$ be a point inside of that ball.
Then any line in $\R P^3$ passing through $x$ intersects the
sphere component of $\R B$ in 2 points.
But it also intersects the plane component of $\R B$
because of homology reasons.
Therefore, $\R B$ is the full inverse image of $\R P^2$
under the projection $\C B\to\C P^2$ from $x$.
The lemma now follows since $\R P^2$ is dual to $w_2(\C P^2)$.
\end{pf}

This lemma allows us to apply the results of \cite{Mi}.
By Theorem 1 and Addendum 1 of \cite{Mi} $\R A$ bounds
such a surface $B_1$ in $\R B$ that
\begin{gather*}
\chi(B_1)\equiv 5\pmod{8},\
\text{if the number of components of $\R A$ is 5,}\\
\chi(B_1)\equiv 4\ \text{or}\ 6\pmod{8},\
\text{if the number of components of $\R A$ is 4.}
\end{gather*}

Therefore if $<\alpha\sqcup 1\!<\!\beta\!>>_{\R P^2}\sqcup<1>_{S^2}$
is realisable then
\begin{gather*}
2+\alpha-\beta\equiv 5\pmod{8},\ \text{if $\alpha+\beta=3$,}\\
2+\alpha-\beta\equiv 4\ \text{or}\ 6\pmod{8},\
\text{if $\alpha+\beta=2$.}
\end{gather*}
Neither of these congruences is possible if $\beta>0$.

To apply the congruence for $\chi(B_1)$ to
$<\alpha\sqcup 1\!<\!\beta\!>>_{\R P^2}\sqcup<\emptyset>_{S^2}$
we have to consider two possibilities
for $B_1$ --- the sphere component of $\R B$ is either contained
in $B_1$ or not.
We get
\begin{gather*}
1+\alpha-\beta\equiv 5\pmod{8}\ \text{or}\
3+\alpha-\beta\equiv 5\pmod{8},\ \text{if $\alpha+\beta=4$,}\\
1+\alpha-\beta\equiv 4\ \text{or}\ 6\ \text{or}\
3+\alpha-\beta\equiv 4\ \text{or}\ 6,\
\text{if $\alpha+\beta=3$.}
\end{gather*}
This is not possible if $\alpha>0$ and $\beta>1$.

\subsection{Constructions}
\begin{figure}
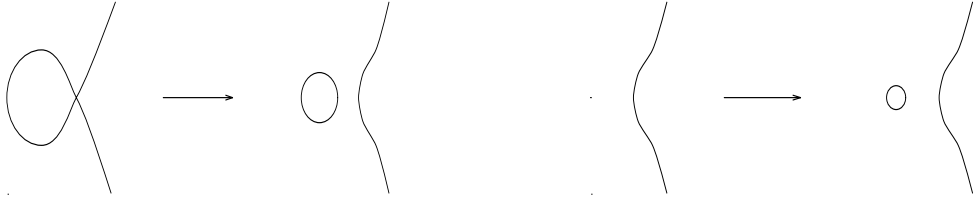

\centerline{\psfig{figure=rhh.eps,height=1in,width=2in}\hspace{1in}
\psfig{figure=rhe.eps,height=1in,width=2in}}
\caption{\label{rh} Construction of \protect{$\R H$}.}
\end{figure}

To construct
$<3\sqcup1\!<\! 1\!>>_{\R P^2}\sqcup<\emptyset>_{S^2}$
we perturb a union of two plane sections of $\R B$.
\begin{lem}
There exists a plane section $\R H$ of $\R B$
which consists of two components which both
belong to the plane component of $\R B$.
\end{lem}
\begin{pf}
Consider a generic plane $P$ in $\R P^3$ tangent to
the plane component of $\R B$.
The intersection $P\cap\R B$ is a singular curve with an ordinary
double point.
Varying $P$ a little we can make $\R H=P\cap\R B$
non-singular and non-connected (see Fig. \ref{rh}).
\end{pf}

\begin{figure}
\centerline{\psfig{figure=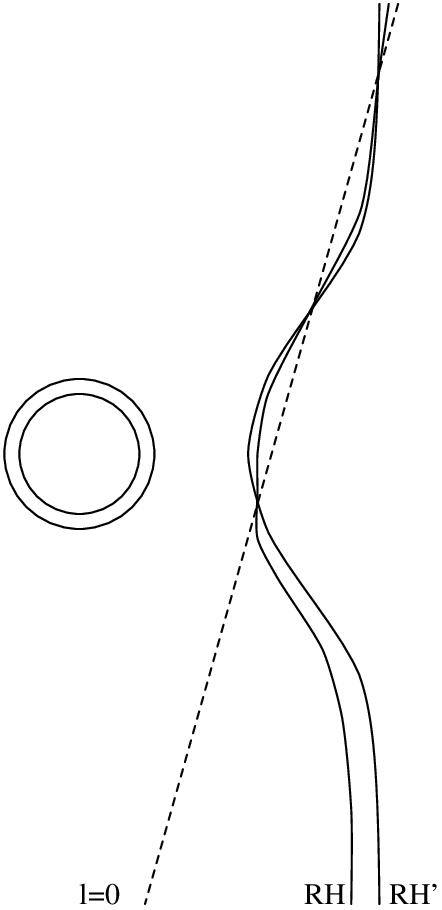,height=2in,width=1.5in}\hspace{0.7in}
\psfig{figure=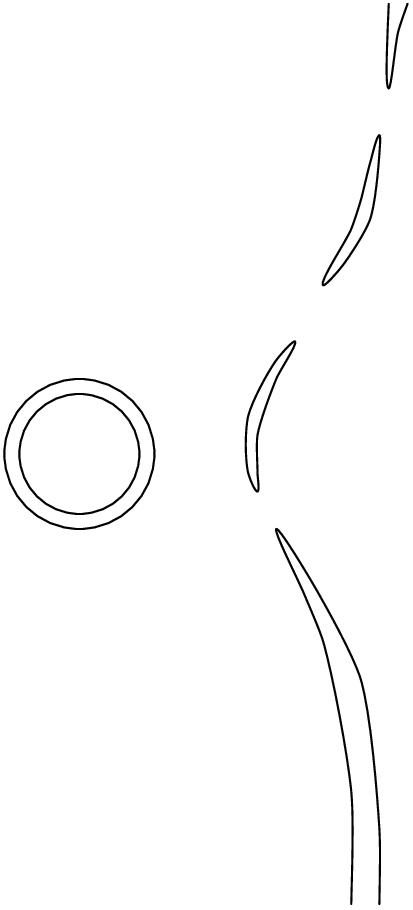,height=2in,width=1.5in}}
\caption{\label{harcub}
\protect{$<3\sqcup1\!<\! 1\!>>_{\R P^2}\sqcup<\emptyset>_{S^2}$}.}
\end{figure}

Let $l=0$ be the equation of the plane in $\R P^3$ intersecting $\R H$
in three distinct points.
Let $\R H'\subset\R B$ be the curve
defined by equation $h+\epsilon l=0$,
where $h$ is the equation of $\R H$ and $\epsilon>0$ is small.
Perturbing $\R H\cup\R H'$ we get the topological arrangement
$<3\sqcup1\!<\! 1\!>>_{\R P^2}\sqcup<\emptyset>_{S^2}$
(see Fig. \ref{harcub}).
The same technique produces
$<2\sqcup 1\!<\!1\!>>_{\R P^2}\sqcup <\emptyset>_{S^2}$
and
$<1\sqcup 1\!<\!1\!>>_{\R P^2}\sqcup <\emptyset>_{S^2}$
if $l=0$ intersects $\R H$ in 2 points (one is a tangent point)
and in 1 point respectively.
This constructs the type a of Theorem 1 with the exception of
$<1\!<\!1\!>>_{\R P^2}\sqcup <\emptyset>_{S^2}$.

\begin{figure}
\centerline{\psfig{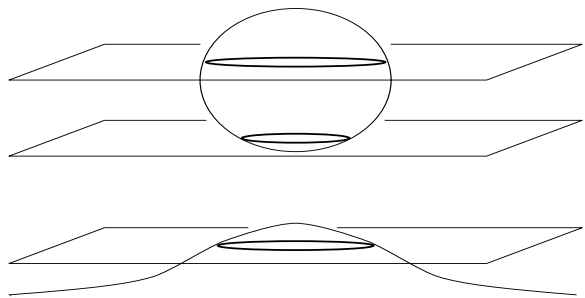}}
\caption{\label{par} Curves on \protect{$x^2+y^2=z(z^2-1)$}.}
\end{figure}

We construct $<1\!<\!1\!>>_{\R P^2}\sqcup <\emptyset>_{S^2}$,
$<1\!<\!1\!<\!1\!>>>_{\R P^2}\sqcup <\emptyset>_{S^2}$ (type d),
$<1\!<\!1\!>>_{\R P^2}\sqcup <1>_{S^2}$ (type e) and
$<\emptyset>_{\R P^2}\sqcup <1\sqcup 1\!<\!1\!>>_{S^2}$ (type f)
on the cubic surface $\R B\subset\R P^3$ given in the affine coordinates
$(x,y,z)$ in $\R^3\subset\R P^3$ by equation
$$x^2+y^2=z(z^2-1).$$
The required topological arrangements are cut on $\R B$ by
equations $(z-C_1)(z-C_2)(z-C_3)=0$ where we vary constants $C_1$,
$C_2$ and $C_3$ (see Fig.\ref{par}).

To construct types b and c
we perturb the union of a plane $P\subset\R P^3$
and an ellipsoid $\R Q\subset\R P^3$
such that $\R H=P\cap\R Q$ is a non-singular non-empty conic.
If $p=0$ is the equation of the plane $P$ and $q=0$
is the equation of the ellipsoid $\R Q$
then $\R B$ is given by equation
$$pq+\epsilon s=0,$$
where $s$ is an auxiliary cubic polynomial which
cuts a non-singular curve $\R S\subset\R Q$
not intersecting $P\cap\R Q$
and $\epsilon\neq 0$ is small.

\begin{figure}
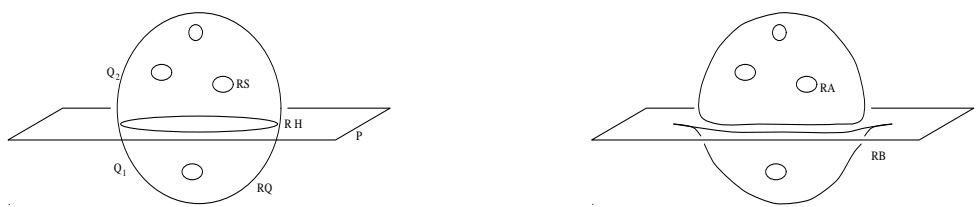

\centerline{\psfig{figure=perturb1.eps,height=1in,width=2in}
\hspace{1in}\psfig{figure=perturb2.eps,height=1in,width=2in}}
\caption{\label{perturb}
Construction of the types b and c of Theorem 1.}
\end{figure}

\begin{figure}
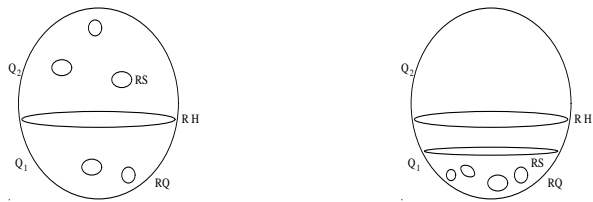

\centerline{\psfig{figure=rs1.eps,height=1in,width=1in}
\hspace{1in}\psfig{figure=rs2.eps,height=1in,width=1in}}
\caption{\label{rs} Auxiliary curves for the types c and b
of Theorem 1.}
\end{figure}

Let $\R A=\R B\cap\R Q$.
The topological type of $(\R B,\R A)$ depends only
on the mutual position of $\R H$ and $\R S$ in $\R Q$
and the sign of $\epsilon$ (see Fig. \ref{perturb}).
The curve $\R H$ splits $\R Q$ into 2 disks $Q_1$
and $Q_2$.
For the arrangements of type c we need
to construct for any $\alpha$, $\beta$, $\alpha+\beta\le 5$,
a curve $\R S$ such that $\R S\cap Q_1$ bounds
$\alpha$ disjoint disks in $Q_1$ and $\R S\cap Q_2$
bounds $\beta$ disjoint disks in $Q_2$.
For the arrangements of type 1.{\bf b} we need
to construct for any $2\le\alpha\le 4$
a curve $\R S$ such that
$\R S\cap Q_1$ has a component which bounds a disk in $Q_1$
which contains all other components of $\R S$ (see Fig. \ref{rs}).

\begin{figure}
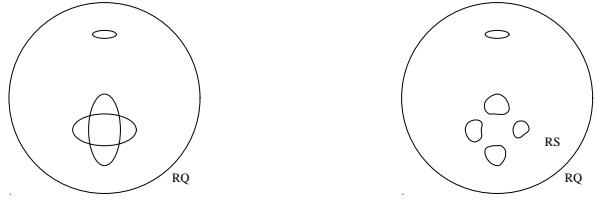

\centerline{\psfig{figure=rsco1.eps,height=1in,width=1in}
\hspace{1in}\psfig{figure=rsco2.eps,height=1in,width=1in}}
\caption{\label{rsco} Construction of \protect{$\R S$}.}
\end{figure}

To construct $\R S\subset\R Q$ consisting of 5 components
we petrurb the union of 3 plane sections of $\R Q$ which
intersect as shown on Fig. \ref{rsco}.
To construct $\R S\subset\R Q$ consisting of 4, 3, 2, 1 or 0
components we just make some of the plane sections (or their
intersections) on Fig. \ref{rsco} imaginary.
It is easy to find all the needed positions of the plane section
$\R H$ directly from Fig. \ref{rsco}.

\section{Proof of Theorem 2}
\begin{figure}
\centerline{\psfig{figure=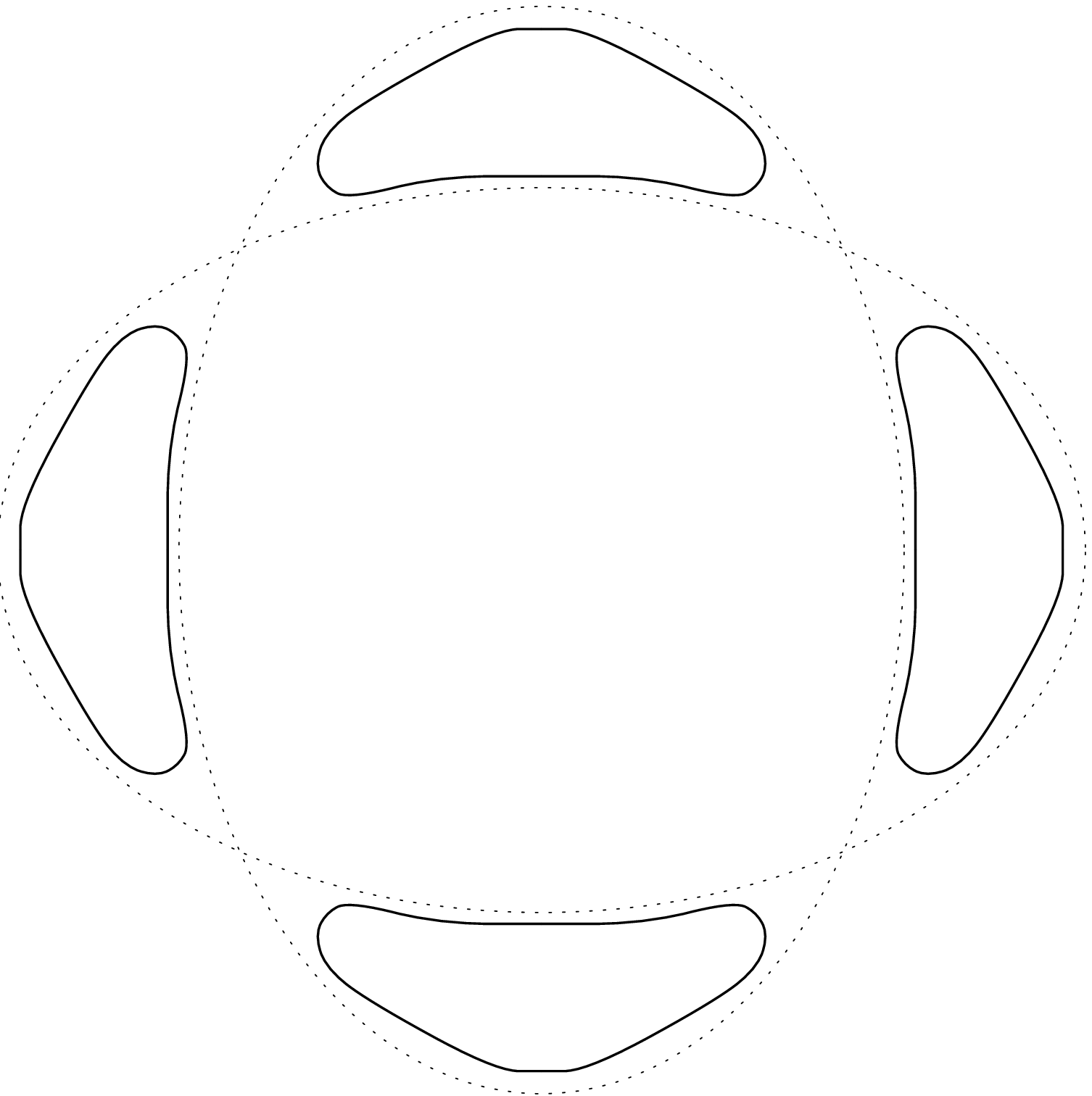,height=1in,width=1in}
\hspace{1in}\psfig{figure=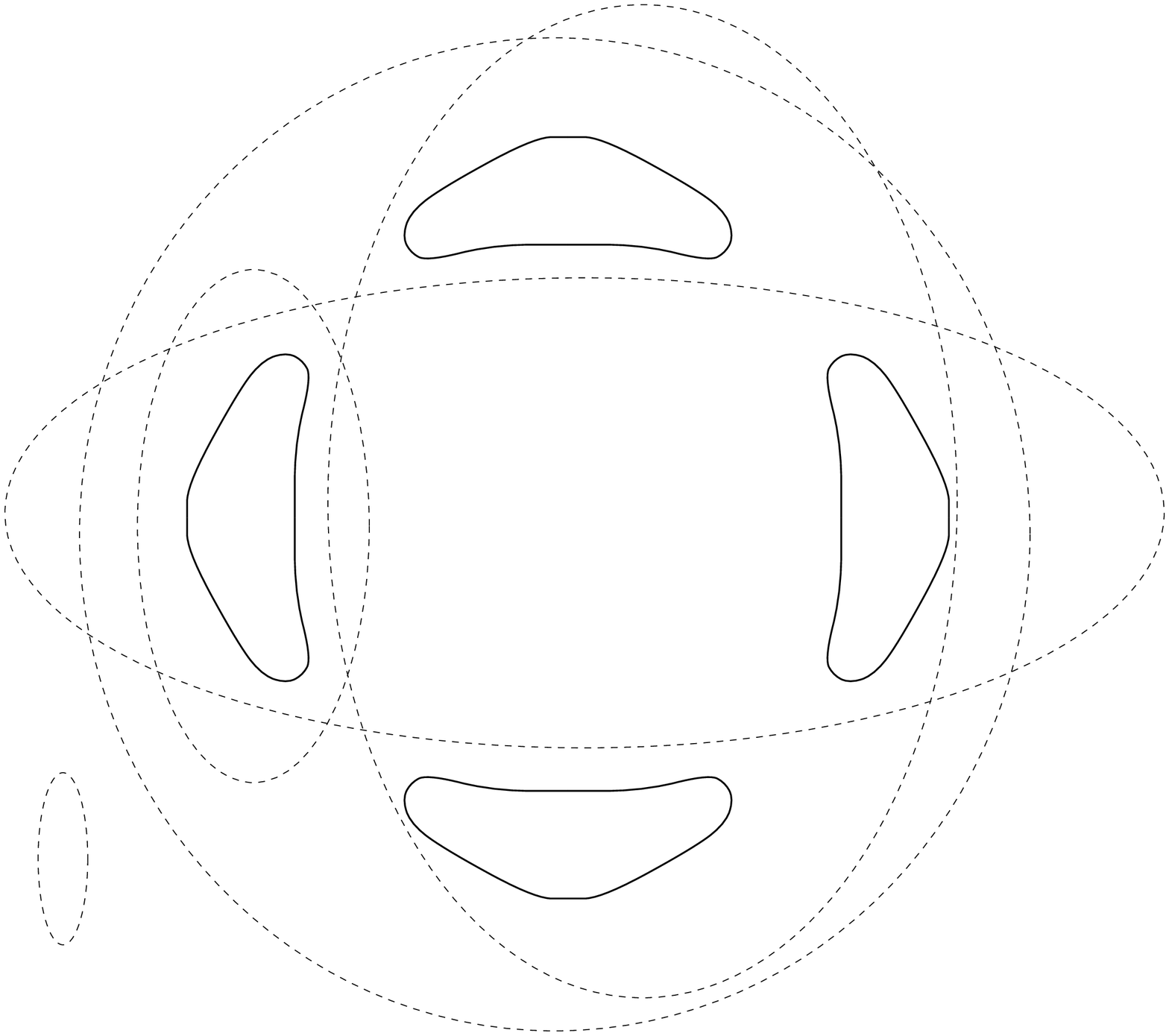,height=1in,width=1in}}
\caption{\label{thm2} Different arrangements of a quartic and a conic.}
\end{figure}

All the needed restrictions for Theorem 2 are
provided by Lemma \ref{arnold}.
To construct the curves we obtain
the cubic surface $\R B$ as the real part
of the blowup $\beta:\C B\to\C P^2$
in 6 imaginary points of $\C P^2$
(these 6 points come in 3 complex conjugate pairs
and they do not belong to the same conic).
Thus to construct a curve $\R A\subset\R B$ of degree 6
it suffices to construct a plane curve $\R C\subset\R P^2$
such that the singular points of its complexifications $\C C$
are 6 imaginary ordinary double points which do not belong
to the same conic.
It is easy to construct a degree 6 curve decomposing to
the union of a quartic and a conic transverse to each other
and with no real intersection points for every topological
arrangement listed in Theorem 2 (see Fig. \ref{thm2}).
This curve has 8 imaginary ordinary double points which
belong to the same conic.
To get $\R C$ we perturb this curve so that 2 double
points disappear and the other 6 survive and move to a
general position (cf. \cite{GU}).

\section{Proof of Theorem 3}
\subsection{Restrictions}
Lemma \ref{arnold} implies that if $B_-$ contains more
than one component of non-positive Euler characteristic
then $(\R B,\R A)$ is of type 8 of Theorem 3.
Otherwise the type of $(\R B,\R A)$ is
\begin{gather*}
<\alpha\sqcup S^2_{\beta+\gamma},\beta\sqcup k\R P^2_{\alpha+\gamma}>,\\
<\alpha\sqcup jT^2_{\beta+\gamma},\beta\sqcup k\R P^2_{\alpha+\gamma}>,\\
\text{or}\
<\alpha\sqcup 2j\R P^2_{\beta+\gamma},\beta\sqcup k\R P^2_{alpha+\gamma}>,\\
\text{where}\ \alpha\ge 0,\ \beta\ge 0,\ \gamma\ge 1,
\alpha+\beta+\gamma\le 5.
\end{gather*}
The number $k$ is determined by the equality
$\chi(B_+)+\chi(B_-)=\chi(\R B)=-5$.
We get $k=9-2j-2\gamma$ (we assume that $j=0$ in the first case).
Since $k\ge 0$ we obtain
$$\gamma<4-j.$$
Therefore to get all the restrictions for Theorem 3
we have to prove that the following 22 types are not realisable.

\begin{itemize}
\item
$<\alpha\sqcup S^2_{\beta+\gamma},
\beta\sqcup (9-2\gamma)\R P^2_{\alpha+\gamma}>$,
where $(\alpha,\beta,\gamma)=(0,3,1)$, $(0,4,1)$,
$(1,2,1)$, $(1,3,1)$, $(2,2,1)$ or $(4,0,1)$.
These types are not realisable because of the Rokhlin (\cite{R})
and Kharlamov-Gudkov-Krakhnov (\cite{Kh}, \cite{GK}) congruences
\begin{gather*}
\chi(B_+)=\alpha-\beta+1\equiv 3\pmod{8},\
\text{if $\alpha+\beta=4$ and $\gamma=1$},\\
\chi(B_+)=\alpha-\beta+1\equiv 2\ \text{or}\ 4\pmod{8},\
\text{if $\alpha+\beta=3$ and $\gamma=1$}.
\end{gather*}
We apply the congruences in
the form of Theorem 7.1 of \cite{Mi}
which is more convenient for our application
(note that it is the condition that $\gamma=1$
which implies the hypothesis $e=0$ of the theorem).

\item
$<\alpha\sqcup T^2_{\beta+\gamma},
\beta\sqcup (7-2\gamma\R P^2_{\alpha+\gamma}>$,
where $(\alpha,\beta,\gamma)=(1,3,1)$
or $(3,1,1)$.
These types are not realizable since Theorem 7.4.d of \cite{Mi}
implies that
$$
\alpha-\beta-1\equiv 3\pmod{4},\
\text{if $\alpha+\beta=4$ and $\gamma=1$.}
$$

\item
$<\alpha\sqcup 2T^2_{\beta+\gamma},
\beta\sqcup (5-2\gamma)\R P^2_{\alpha+\gamma}>$,
where $(\alpha,\beta,\gamma)=(0,4,1)$, $(2,2,1)$
or $(4,0,1)$.
These types are not realizable since Theorem 7.4.d of \cite{Mi}
implies that
$$
\alpha-\beta-3\equiv 3\pmod{4},\
\text{if $\alpha+\beta=4$ and $\gamma=1$.}
$$

\item
$<\alpha\sqcup 3T^2_{\beta+1},
\beta\sqcup\R P^2_{\alpha+1}>$,
where $(\alpha,\beta)=(4,0)$, $(1,1)$, $(2,1)$, $(3,1)$,
$(1,2)$, $(2,2)$ or $(1,3)$.
These types are not realizable since Theorem 1
of \cite{Mi} implies
$$
\chi(B_-)\equiv 3+B\pmod{8},
$$
where $B$ is the Brown invariant of the Guillou-Marin form
of $\C A/\conj\cup B_-$.
But $B\equiv 1\pmod{8}$ if $\R A$ is an M-curve
and $B\equiv 0\ \text{or}\ 2\pmod{8}$ if $\R A$ is an (M-1)-curve
(cf. Addenda 1.a, 1.b of \cite{Mi})
since the only class in
$H_1(B_-;\Z_2)$ which does not vanish in $H_1(\R B;\Z_2)$ is the
class dual to $w_1(\R B)$ and the value of the Guillou-Marin form of
$\C A/\conj\cup B_-$ on such a class is 1.
Therefore
\begin{gather*}
\beta-\alpha\equiv 4\pmod{8},\
\text{if $\alpha+\beta=4$ and $\gamma=1$},\\
\beta-\alpha\equiv 3\ \text{or}\ 5\pmod{8},\
\text{if $\alpha+\beta=3$ and $\gamma=1$.}
\end{gather*}
This rules out the pairs
$(\alpha,\beta)=(2,1),(3,1),(1,2),(2,2),(1,3)$.
In a similar way Addendum 1.c of \cite{Mi} implies
that $(\alpha,\beta)=(1,1)$ is not realizable by
a curve of type II (i.e. if $\C A/\conj$ is non-orientable)
since $\beta-\alpha\equiv 0\pmod{8}$ in this case.

To see that $(\alpha,\beta)=(4,0)$
is not realizable
we apply a version of the Rokhlin complex orientation
formula \cite{R1}.
Note that $B_+$ in this case is a disjoint union of
4 disks and a punctured sphere with 3 handles
and therefore $F=\C A/\conj\cup B_+\subset\C B/\conj\approx S^4$
(see \cite{L}) is orientable.
Therefore $F.F=0$ but on the other hand
$F.F=\frac{1}{2}\C B.\C B-2\chi(B_+)=6+2=8$.

A similar argument shows that $(\alpha,\beta)=(1,1)$
is not realizable by a curve of type I
(i.e. if $\C A/\conj$ is orientable).
In this case $B_+$ is a disjoint union of a disk and
a sphere with 3 handles punctured two times and
$B_-$ is a disjoint union of a disk and a projective plane
punctured two times.
Denote the disk component of $B_-$ by $D_-$
the disk component of $B_+$ by $D_+$ and
the union of $D$ and the non-disk component of $B_+$
by $E$.
Note that $G=\C A/\conj\cup D_+\cup E\cup D_-$ is
a $\Z$-cycle in $\C B/\conj\approx S^4$ and thus $G.G=0$.
If the orientations of $\R A=\dd\C A/\conj=\dd B_+$
induced from $\C A/\conj$ and from $B_+$ agree then
$G.G=\frac{1}{2}\C B.\C B-2(\chi(D_+)+\chi(E)+\chi(D_-))+4=16\neq 0$.
If the orientations do not agree then
$G.G=\frac{1}{2}\C B.\C B-2(\chi(D_+)+\chi(E)+\chi(D_-))-4=8\neq 0$.

\item
$<2\sqcup 2\R P^2_3,2\sqcup 5\R P^2_3>$.
This type is not realisable because of Theorem 7.4.c of \cite{Mi},
since $\chi(B_+)=-1$ and $B_+$ is not orientable.

\item
$<\alpha\sqcup 6\R P^2_{\beta+1},\beta\sqcup\R P^2_{\alpha+1}>$,
where $(\alpha,\beta)=(2,1)$, $(3,1)$ or $(2,2)$.
The pair $(2,1)$ is ruled out by Theorem 7.4.c of \cite{Mi}
applied to $B_-$ since $\chi(B_-)=\beta-\alpha=-1$ and
$B_-$ is not orientable.
The pairs $(3,1)$ and $(2,2)$ are ruled out by Theorem 7.4.b
of \cite{Mi}.
\end{itemize}

\subsection{Constructions}
\label{thm3con}
\begin{figure}
\centerline{\psfig{figure=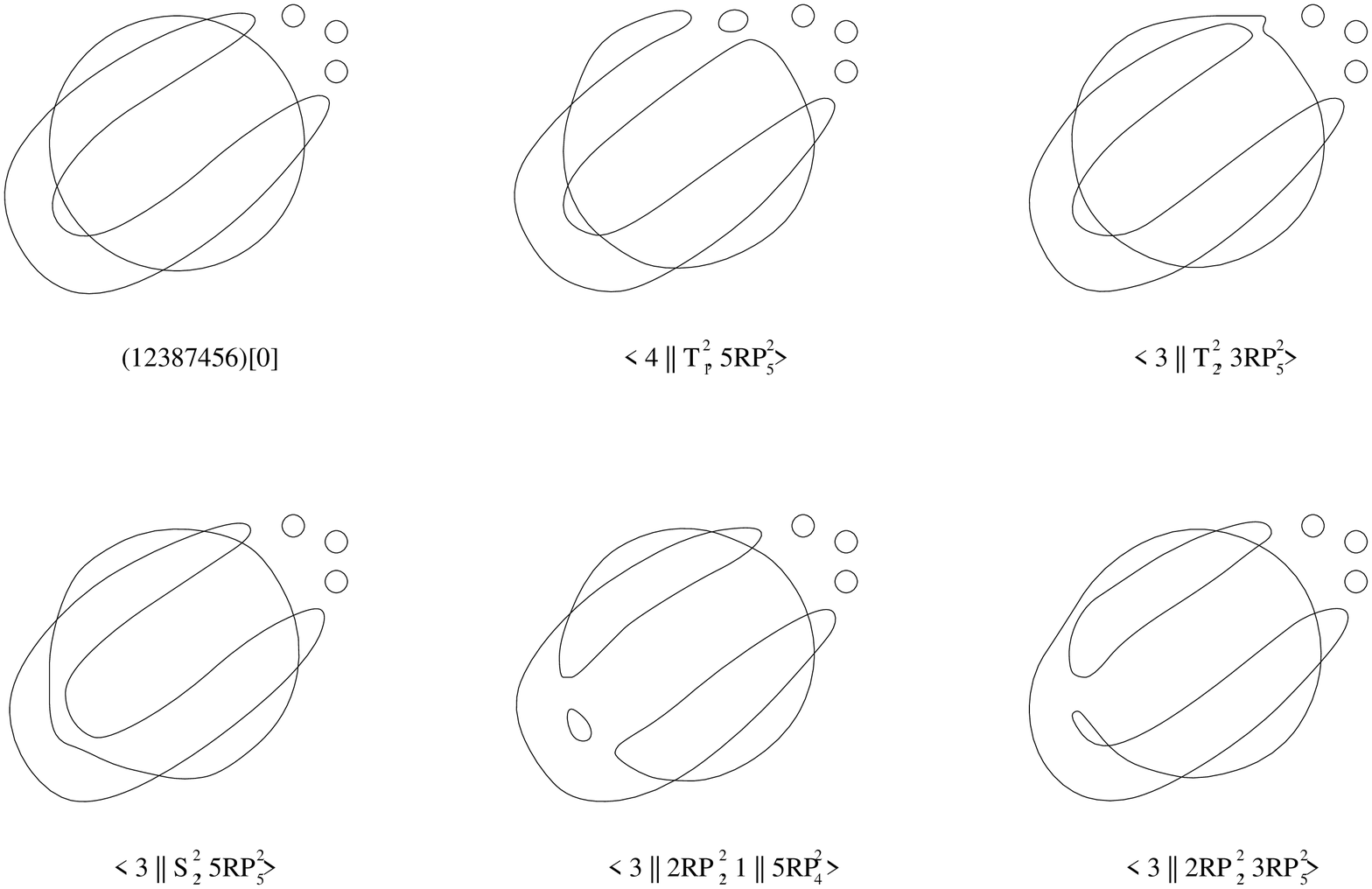,height=3.4in,width=5in}}
\caption{\label{beispiel} Curves produced by $(12387456)[0]$}
\end{figure}

Similarily to the proof of Theorem 2
we construct the curves for Theorem 3 perturbing the union
of a quartic and a conic.
This time we have to use a quartic and a conic which intersect
in at least 6 real points since we need 6 points in $\R P^2$
to blow up to get the cubic surface diffeomorphic to $7\R P^2$.
We use Polotovskii's catalog \cite{P} which lists
arrangements of a conic and a quartic
which intersect in 8 real points.
Then we smooth 2 out of the 8 points so that the other 6
move to a general position with respect to each other (cf. \cite{GU}).
One curve in Polotovskii's catalog may produce more
than one type for Theorem 3
because of the ambiguity in the choice of the 2 points to smooth
and in (real) smoothing of these points.
\begin{ex}
The curve encoded by $(12387456)[0]$ in Polotovskii's catalog \cite{P}
produces $<3\sqcup S^2_2,5\R P^2_5>$, $<4\sqcup T^2_1,5\R P^2_5>$,
$<3\sqcup T^2_2,3\R P^2_5>$, $<3\sqcup 2\R P^2_2,1\sqcup 5\R P^2_4>$
and $<3\sqcup 2\R P^2_2,3\R P^2_5>$ out of M-curves (see Fig. \ref{beispiel})
and $<3\sqcup T^2_1,5\R P^2_4>$ and $<3\sqcup 2\R P^2_1,5\R P^2_4>$
out of $(M-1)$-curves.
\end{ex}

\begin{figure}[b]
\centerline{\psfig{figure=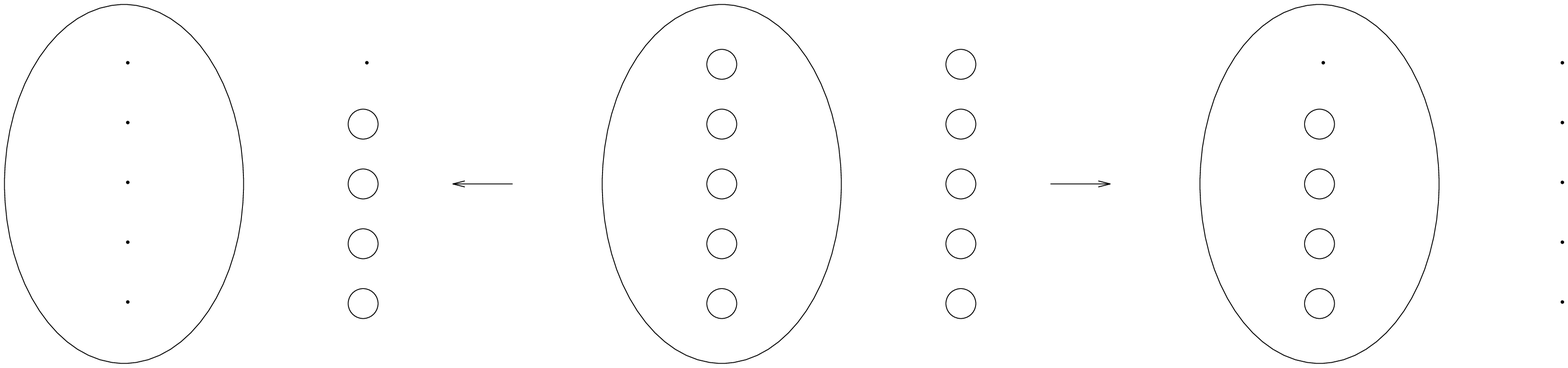,height=1in,width=4.5in}}
\caption{\label{gud}
$<2\sqcup\R P^2_5,4\sqcup 5\R P^2_1>$
and $<4\sqcup 6\R P^2_1,\R P^2_5>$.}
\end{figure}

\begin{figure}
\centerline{\psfig{figure=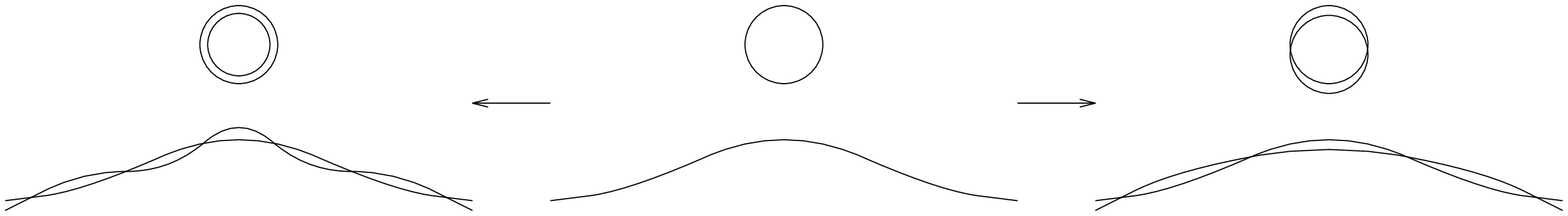,height=0.7in,width=5in}}
\caption{\label{arn}
$<1\sqcup 3T^2_2,S^2_2\sqcup\R P^2_1>$
and $<2T^2_3,S^2_2\sqcup\R P^2_1>$.}
\end{figure}

The catalog \cite{P} lists only arrangements of M-decomposing curves
thus the quartic in the union consists of 4 ovals.
But if an oval of the quartic is disjoint from the conic then
we can remove it so that the resulting arrangement is still
realizable by a decomposing curve of degree 6.
This means that if a triple $(\alpha,\beta,\gamma)$ is obtained
from a decomposing curve then $(\alpha',\beta',\gamma)$ with
$\alpha'<\alpha$ and $\beta'<\beta$ can be also obtained in this way.
In the following tables
we indicate the code of a union of a conic and a quartic in \cite{P}
for arrangements in Theorem 3 with maximal $(\alpha,\beta)$.
\vspace{10pt}\newline
{\bf 1.} $<\alpha\sqcup S^2_{\beta+\gamma},
\beta\sqcup(9-2\gamma)\R P^2_{\alpha+\gamma}>$,\\
\begin{tabular}{|r|r|r|r|r|}
\hline
\ & $\beta=0$ & $\beta=1$ & $\beta=2$ & $\beta=3$\\
\hline
\ & \ & $\gamma=4$ & $\gamma=3$ & $\gamma=2$ \\
$\alpha=0$ & \ & (12)(34)(56)(78) &
(1867)(3452)[2] & (18276543)[3]\\
\ & \ & \ & $\gamma=1$ & \ \\
\ & \ & \ & * & \ \\
\hline
\ & $\gamma=4$ & $\gamma=3$ & $\gamma=2$ & \ \\
$\alpha=1$ & (12)(34)(56)(78) & (1678)(2345)[1] &
(18723456)[2] & \ \\
\hline
\ & $\gamma=3$ & $\gamma=2$ & \ & \ \\
$\alpha=2$ & (145678)(23)[0] & (12345678)[1] & \ & \ \\
\hline
\ & $\gamma=2$ & $\gamma=1$ & \ & \ \\
$\alpha=3$ & (12387456)[0] & (12345678)[1] & \ & \ \\
\hline
\end{tabular}\vspace{10pt}\\
{\bf 2.} $<\alpha\sqcup T^2_{\beta+\gamma},
\beta\sqcup(7-2\gamma)\R P^2_{\alpha+\gamma}>,$\\
\begin{tabular}{|r|r|r|r|r|r|}
\hline
\ & $\beta=0$ & $\beta=1$ & $\beta=2$ & $\beta=3$ & $\beta=4$ \\
\hline
\ & \ & \ & $\gamma=3$ & $\gamma=2$ & $\gamma=1$ \\
$\alpha=0$ & \ & \ & (1876)(2345)[2] &
(18743256)[3] & (18276543)[3] \\
\hline
\ & \ & $\gamma=3$ & $\gamma=2$ & \ & \ \\
$\alpha=1$ & \ & (187654)(23)[1] & (18723456)[2] & \ & \ \\
\hline
\ & $\gamma=3$ & $\gamma=2$ & $\gamma=1$ & \ & \ \\
$\alpha=2$ & (1867)(3452)[2] & (18765234)[1] & (18723456)[2] & \ & \ \\
\hline
\ & $\gamma=2$ & \ & \ & \ & \ \\
$\alpha=3$ & (12387456)[0] & \ & \ & \ & \ \\
\hline
\ & $\gamma=1$ & \ & \ & \ & \ \\
$\alpha=4$ & (12387456)[0] & \ & \ & \ & \ \\
\hline
\end{tabular}
\vspace{10pt}\\
{\bf 3.} $<\alpha\sqcup 2T^2_{\beta+\gamma},
\beta\sqcup(5-2\gamma)\R P^2_{\alpha+\gamma}>$,\\
\begin{tabular}{|r|r|r|r|r|}
\hline
\ & $\beta=0$ & $\beta=1$ & $\beta=2$ & $\beta=3$ \\
\hline
\ & \ & \ & \ & $\gamma=2$ \\
$\alpha=0$ & \ & \ & \ & (18743256)[3] \\
\hline
\ & \ & \ & $\gamma=2$ & $\gamma=1$ \\
$\alpha=1$ & \ & \ & (18765432)[2] & (18234765)[2] \\
\hline
\ & \ & $\gamma=2$ & \ & \ \\
$\alpha=2$ & \ & (18765234)[1] & \ & \ \\
\hline
\ & $\gamma=2$ & $\gamma=1$ & \ & \ \\
$\alpha=3$ & (18276543)[3] & (18765234)[1] & \ & \ \\
\hline
\end{tabular}\vspace{10pt}\\
{\bf 4.} $<\alpha\sqcup 3T^2_{\beta+1},
\beta\sqcup\R P^2_{\alpha+1}>$,\\
\begin{tabular}{|r|r|}
\hline
$\beta=0$ & $\beta=4$\\
\hline
$\alpha=3$ & $\alpha=0$\\
(18276543)[3] & (18234567)[3] \\
\hline
\end{tabular}\vspace{10pt}\\
{\bf 5.} $<\alpha\sqcup 2\R P^2_{\beta+\gamma},
\beta\sqcup (7-2\gamma)\R P^2_{\alpha+\gamma}>$,\\
\begin{tabular}{|r|r|r|r|r|r|}
\hline
\ & $\beta=0$ & $\beta=1$ & $\beta=2$ & $\beta=3$ & $\beta=4$ \\
\hline
\ & \ & \ & $\gamma=3$ & $\gamma=2$ & $\gamma=1$ \\
$\alpha=0$ & \ & \ & (1867)(3452)[2] &
(18276543)[3] & * \\
\hline
\ & \ & $\gamma=3$ & $\gamma=2$ & $\gamma=1$ & \ \\
$\alpha=1$ & \ & (1845)(23)(67)[0] & (18723456)[2] &
(18723456)[2] & \ \\
\hline
\ & $\gamma=3$ & $\gamma=2$ & \ & \ & \ \\
$\alpha=2$ & (1845)(3672)[0] & (18432765)[1] & \ & \ & \ \\
\hline
\ & $\gamma=2$ & $\gamma=1$ & \ & \ & \ \\
$\alpha=3$ & (12387456)[0] & (12387456)[0] & \ & \ & \ \\
\hline
\ & $\gamma=1$ & \ & \ & \ & \ \\
$\alpha=4$ & (16254378)[0] & \ & \ & \ & \ \\
\hline
\end{tabular}\vspace{10pt}\\
{\bf 6.} $<\alpha\sqcup 4\R P^2_{\beta+\gamma},
\beta\sqcup (5-2\gamma)\R P^2_{\alpha+\gamma}>$,\\
\begin{tabular}{|r|r|r|r|r|r|}
\hline
\ & $\beta=0$ & $\beta=1$ & $\beta=2$ & $\beta=3$ & $\beta=4$ \\
\hline
\ & \ & \ & \ & $\gamma=2$ & $\gamma=1$ \\
$\alpha=0$ & \ & \ & \ &
(18437625)[3] & (18437625)[3] \\
\hline
\ & \ & \ & $\gamma=2$ & $\gamma=1$ & \ \\
$\alpha=1$ & \ & \ & (18432765)[1] &
(18437625)[3] & \ \\
\hline
\ & \ & $\gamma=2$ & $\gamma=1$ & \ & \ \\
$\alpha=2$ & \ & (18765234)[1] & (18765234)[1] & \ & \ \\
\hline
\ & $\gamma=2$ & $\gamma=1$ & \ & \ & \ \\
$\alpha=3$ & (18276345)[0] & (18276345)[0] & \ & \ & \ \\
\hline
\ & $\gamma=1$ & \ & \ & \ & \ \\
$\alpha=4$ & (18276345)[0] & \ & \ & \ & \ \\
\hline
\end{tabular}\vspace{10pt}\\
{\bf 7.} $<\alpha\sqcup 6\R P^2_{\beta+1},
\beta\sqcup \R P^2_{\alpha+1}>$,\\
\begin{tabular}{|r|r|r|}
\hline
$\beta=0$ & $\beta=3$ & $\beta=4$\\
\hline
$\alpha=4$ & $\alpha=1$ & $\alpha=0$\\
** & (18765432)[2] & (18234567)[3] \\
\hline
\end{tabular}
\vspace{10pt}\\
To construct the arrangements marked by * in the table
we start from a plane curve of degree 6 and collapse
6 of its empty ovals to points using \cite{I}.
To construct $<2\R P^2_5,4\sqcup 5\R P^2_1>$ and
$<4\sqcup 6\R P^2_1,\R P^2_5>$ we use the Gudkov curve
$<5\sqcup 1\!<\!5\!>>$ constructed in \cite{G1}.
To get $<2\R P^2_5,4\sqcup 5\R P^2_1>$
we collapse the 5 inner ovals and one
of the outer ovals;
to construct $<4\sqcup 6\R P^2_1,\R P^2_5>$
we collapse the 5 outer ovals and one
of the inner ovals (see Fig. \ref{gud}).
To construct $<S^2_3,2\sqcup 7\R P^2_1>$
we collapse the 6 outer ovals
of $<6\sqcup 1\!<\!2\!>>$.
\vspace{10pt}\\
{\bf 8.}
$<1\sqcup 3T^2_2,S^2_2\sqcup\R P^2_1>$,
$<2T^2_3,S^2_2\sqcup\R P^2_1>$\vspace{3pt}\\
To construct the plane curves of degree 6 with
6 double points for these arrangements we start from the union of a
two-component cubic in $\R P^2$ and its "parallel" copy which
intersect in 9 disjoint point and smooth 3 out of the 9 nodes of the
union (see Fig. \ref{arn}).

\section{Proof of Theorem 4 and Theorem 5}
\subsection{Restrictions}
Lemma \ref{arnold} implies that unless $(\R B,\R A)$
is of type 6 of Theorem 4 or of type 4 of Theorem 5
it is of type
\begin{gather*}
<\alpha\sqcup S^2_{\beta+\gamma},\beta\sqcup k\R P^2_{\alpha+\gamma}>,\\
<\alpha\sqcup jT^2_{\beta+\gamma},\beta\sqcup k\R P^2_{\alpha+\gamma}>,\\
\text{or}\
<\alpha\sqcup 2j\R P^2_{\beta+\gamma},\beta\sqcup k\R P^2_{\alpha+\gamma}>,\\
\text{where}\ \alpha\ge 0,\ \beta\ge 0,\ \gamma\ge 1,
\alpha+\beta+\gamma\le 5.
\end{gather*}
In the case of Theorem 4 $k=7-2j-2\gamma\ge 0$
and in the case of Theorem 5 $k=5-2j-2\gamma\ge 0$
since $\chi(B_+)+\chi(B_-)=\chi(\R B)$.
To finish the restrictions for Theorem 4 and Theorem 5
we have to prove that $<4\sqcup 2T^2_1,\R P^2_5>$ is not realizable
(for Theorem 4).
Note that $F=\C A/\conj\cup B_+\subset\C B/\conj\approx\overline{\C P}^2$
(see \cite{L})
is orientable and therefore $F.F\le 0$ ($\overline{\C P}^2$ is negative
definite).
On the other hand $F.F=\frac{1}{2}\C B.\C B-2\chi(B_+)=6-2=4$.

\subsection{Constructions}
\begin{figure}
\centerline{
\psfig{figure=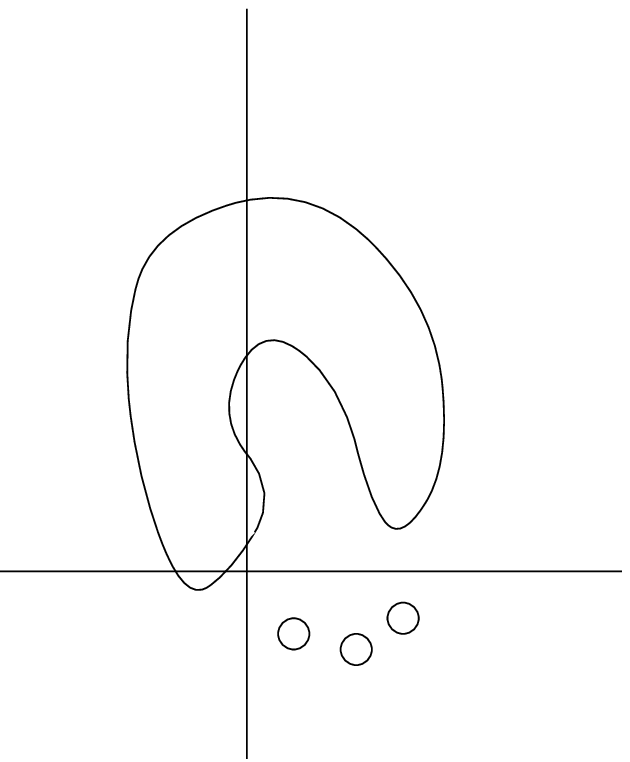,height=1.2in,width=1in}\hspace{0.2in}
\psfig{figure=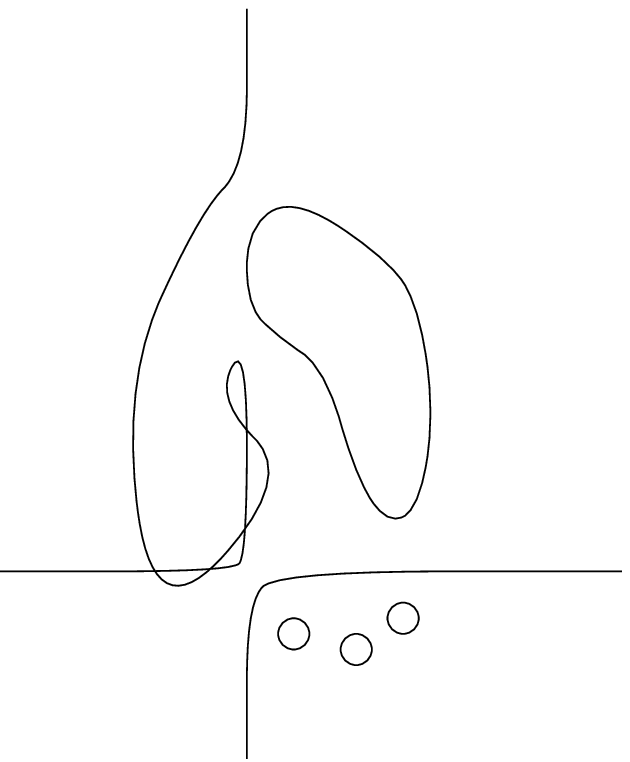,height=1.2in,width=1in}\hspace{0.6in}
\psfig{figure=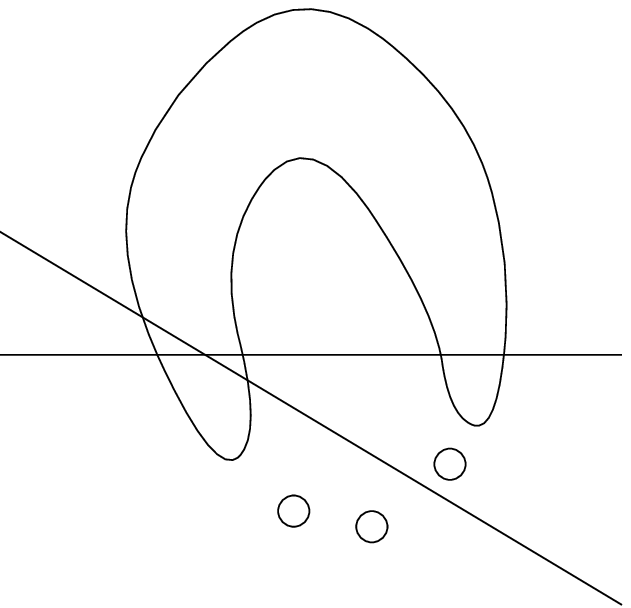,height=1.2in,width=1in}\hspace{0.2in}
\psfig{figure=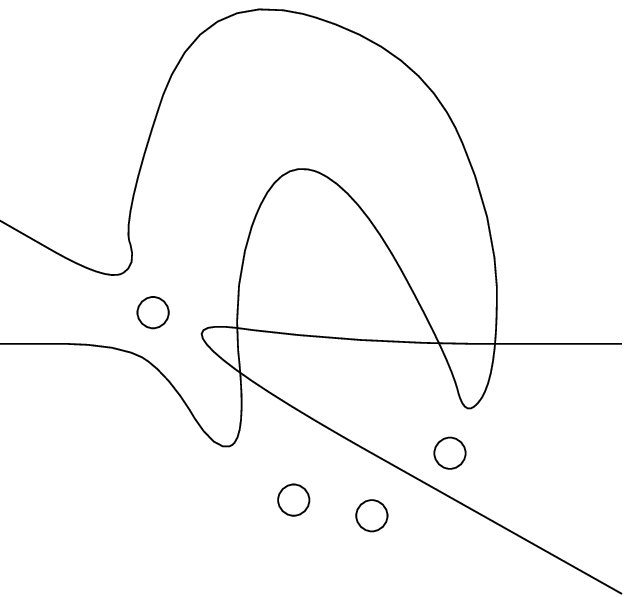,height=1.2in,width=1in}}
\caption{\label{arr}
$<1\sqcup S^2_4,3\sqcup 5\R P^2_2>$ and
$<2\sqcup 2T^2_3,2\sqcup\R P^2_3>$.}
\end{figure}

If $<j\sqcup F\#T^2,n\sqcup G>$ or $<j\sqcup F,n\sqcup G\#T^2>$
is constructed in \ref{thm3con} out of Polotovskii's catalog
then we can construct $<j\sqcup F,n\sqcup G>$ out of
a similar pair of a conic and a quartic
which intersect in 6 real points.
Considering similar pairs of a conic and a quartic
which intersect in 4 real points produces the arrangements
for Theorem 5.
A similar modification works for
$<1\sqcup 2T^2_2,S^2_2\sqcup\R P^2_1>$ and $<T^2_3,S^2_2\sqcup\R P^2_1>$.
This produces all the arrangements for Theorem 4 except for
$<1\sqcup S^2_4,3\sqcup 5\R P^2_2>$ and $<2\sqcup 2T^2_3,2\sqcup\R P^2_3>$.
To construct them we perturb a union of a quartic and two lines
as shown on Fig. \ref{arr}.

\end{document}